\documentclass[12pt]{article}

\begin{document}
\bigskip

\begin{center}

{\Large \bf
 Stochastic Properties of Dynamical Systems Arising from (quantum) Spaces and Actions of  Groups
 }
\end{center}
\smallskip
\begin{center}
{\bf Nikolaj M. Glazunov} \end{center}
\smallskip
\begin{center}
{\rm  National Aviation University, Kiev. }
\end{center} 
\begin{center}
{\rm  Email:} {\it glanm@yahoo.com }
\end{center} \smallskip

{\bf  Abstract:} 

 We review novel results and investigate actions and transformations of  groups and semigroups on 
(quantum) spaces, present  dynamical systems  and   zeta functions arising from these spaces,  actions 
and transformations, discuss their stochastic properties.

\smallskip
{\bf  Keywords:} Dynamical System; Ergodic Transformation; Group Action; Equidistribution; Zeta function; 
Arithmetic Surface.

\section{Introduction}
 
A history of a semigroup and a group action on tori and projective spaces can be found among  other in the
book by A.G. Postnikov \cite{Postnikov2005}, in the paper by I.Ya. Gol'dsheid, G.A. Margulis \cite{Gol'dsheidMargulis}  
and in the supplement by B.M. Gurevich, Ya. G. Sinai \cite{GurevichSinai}   to the Russian translation of the English 
edition of the book by P. Billingsley \cite{Billingsley}.
\par
Dynamical systems have several important complexity measures among which are measure-theoretic, topological and metric 
entropies. A. Kolmogorov defines the complexity of a measure-preserving transformation by generators. A generator for 
a measure-preserving transformation $T$ is a partition $\xi$ with finite entropy such that the set of finite entropy 
partitions subordinate to some $\bigvee_{i=-n}^{n} T^{-i} (\xi)$ is dense in the set of finite entropy partitions 
endowed with the Rokhlin metric. For smooth dynamical systems topological entropy characterizes the total exponential 
complexity of the orbit structure.  Metric entropy with respect to an invariant measure codes the exponential growth rate 
of the statistically significant orbits. In symbolic dynamics complexity of a dynamical system is measured with respect to 
a coding of its orbits.  Dynamical systems can be generated among  other by iteration of maps, $\beta$-transformations, 
Hasse-Kloosterman maps, partitions, group actions.
Here we review novel results and investigate actions and transformations of  groups and semigroups on 
(quantum) spaces, present  dynamical systems  and   zeta functions arising from these spaces,  actions and transformations, 
discuss their stochastic properties.

\section{Dynamical systems from spaces}
It is well known that one-dimensional projective space ${\bf P}^{1}({\bf Q})$ parametrize the set of dynamical systems in 
such a way that for any rational point $Q \in {\bf P}^{1}({\bf Q}), Q = (\frac{a}{b}, 1), a, b \in {\bf Z}, (a,b) = 1$  we 
naturally assiciate dynamical system $({\bf T}, T_Q)$. Here 
${\bf T} = {\bf R}/{\bf Z}, {\bf T}^{\bf Z} = (...,x_{-1}, x_{0}, x_{1},...), x_{i} \in {\bf T}, X = \{{\bf x}=(x_{k}): 
bx_{k + 1} = ax_{k} \; \mbox{for all}\; k \in {\bf Z}\},  T_Q : X \to X$. More generally, for any primitive polynomial 
$g(x) \in {\bf Z}[x]$ of degree $d \ge 1$ it is possible to construct its Frobenius and companion matrices and define a 
homeomorphism $T_F$ of a compact $d-$dimensional subgroup of ${\bf T}^d$.
These considerations can be extended to elliptic curves \cite{DEMW:DS} and to abelian varieties.
For elliptic curves authors of the paper \cite{DEMW:DS} implement these by the following way. Let $q \in {\bf Q}_p$ and 
$\log^{+}x$ denotes $\max\{\log x,0\}.$ For a generic element $x$ of ${\bf Z}_p$ authors define $q$-transformation
$T_{q}(x)$ (a $p$-adic analogue of the $\beta$-transformation). Then the topological entropy of the $p$-adic 
$\beta$-transformation is given by $h(T_{q}) = \log^{+}|q|_p 
\;$ (\cite{DEMW:DS}, Theorem 4.1). If $|q|_p \geq 1$ then the map $T_{q}$ is
ergodic with respect to Haar measure for $|q|_p > 1$ and is
not ergodic for $|q|_p = 1$ (\cite{DEMW:DS}, Theorem 4.2). Let $Per_{n}(T_{q})$
denotes the subgroup of ${\bf Z}_p$ consisting of elements
of period $n$ under $T_{q}$. Let $U$ be the set of unit roots
of ${\bf Q}_p$ and $q \in {\bf Q}_p \setminus U$. Then
$$ \log |Per_{n}(T_{q})| = n \log^{+} |q|_{p}. $$
(\cite{DEMW:DS}, Theorem 4.3). The authors use the topological entropy and 
measure theoretical arguments based on volume growth rate and
arithmetic of ${\bf Z}_p.$ \\ 
Let $Q$ be a rational point of an elliptic curve over ${\bf Q}$ 
and let $\hat{h}(Q)$
be the global canonical height on rational points of the
elliptic curve. Then with the definitions and assumptions of the 
paper \cite{DEMW:DS} and $q = a/b = x(Q)$, (i) the entropy of 
$T_Q$ is given by $h(T_{Q}) = 2\hat{h}(Q),$ and (ii) the asymptotic
growth rate of the periodic points is given by the division
polynomial $\nu_{n}(x) $:
$ \log |Per_{n}(T_{Q}| \sim \log |b^{n} \nu_{n}(q)|$ as 
$n \rightarrow \infty.$ (\cite{DEMW:DS}, Theorem 5.2). In the case authors use
also the elliptic analogue of Baker's theorem, which described
in paper \cite{D:DS}  and in paper \cite{EW:EM} .

\section{Dynamical systems on probability spaces}
Let $(X, B, \mu, T)$ be a dynamical system on standard probability space with $T : X \to X$ is measurable, almost
 surely one to one, preserves $\mu$, for which it is an ergodic transformation.
 Random dynamical systems relate a partial case of  bundle dynamical systems by I. Cornfeld, S.  Fomin, and 
Ya.  Sinai \cite{CFS:ES}. Measurable partition of the space $X$ transforms the initial random dynamical system into 
a symbolic dynamical system. Below we will present novel symbolic dynamical systems and their applications.

\section{Rigid and weakly mixing ergodic transformations}
In papers \cite{AS:GT}  and \cite{JKLS:PA}  authors present resent results on genericity of rigid and multiply recurrent 
infinite measure preserving and nonsingular transformations  and on measurable sensitivity. In the paper \cite{JKLSS:ES}
authors investigate properties of uniformly rigid transformations and analyze the compatibility of uniform rigidity and 
measurable weak mixing along with some of their asymptotic convergence properties.
All spaces of the paper under review are considered simultaneously as topological spaces and as measure spaces. Presented 
results concern either the measurable dynamics on the spaces or the interplay between the measurable and topological dynamics. 
The notion of uniform rigidity was introduced as a topological version of rigidity by S. Glasner and D. Maon \cite{GM:ET}.
Authors of the paper \cite{JKLSS:ES} considers functional analytic properties of uniform rigidity that is similar to the properties 
of rigidity.
Theorem 1 (\cite{JKLSS:ES}). Every totally ergodic finite measure-preserving transformation on a Lebesgue space has a 
representation that is not uniformly rigid, except in the case where the space consists of a single atom.

The proof of the theorem connects with results of authors of the theorem that uniform rigidity and weak mixing are mutually 
exclusive notions on a Cantor set, and follows from the Jewett-Krieger Theorem by K. Peterson \cite{Pe:ETh}.

\section{Superrigidity for groups}
The concept of superrigidity was introduced by G. D.  Mostow \cite{Mo:SR}
and by G. A. Margulis \cite{Ma:DS}
in the context of studying the structure of lattices in rank one and higher rank Lie groups respectively.
The notion of property (T) for locally compact groups was defined by 
D. Kazhdan \cite{Ka:FA}  and the notion of relative property (T)  
for inclusion of countable groups $\Gamma_0 \subset \Gamma$ was defined by 
G. Margulis \cite{Ma:ES}. Now considerthe orbit equivalence (OE) superrigidity. One of 
 the first result of this type of superrigidity was obtained by A. Furman \cite{Fu:ES},
 who combined the cocycle superrigidity by 
R. Zimmer \cite{Zi:ET}
with ideas from geometric group theory to show that the actions 
$SL_{n}({\bf Z})$ on  ${\bf T}^n (n \ge 3)$ are OE superrigid. The deformable actions of rigid 
groups are OE superrigid by S. Popa \cite{Po:IC}.  
The main result of the paper by A. Ioana \cite{Io:CS} is the Theorem A on orbit equivalence (OE) superrigidity. As a consequence 
of Theorem A the author of the paper \cite{Io:CS} can constructs uncountable many non-OE profinite actions for the arithmetic groups 
$SL_n ({\bf Z}) (n \ge 3)$, as well as for their finite subgroups, and for the groups 
$SL_m ({\bf Z}) \times {{\bf Z}}^m (m \ge 2)$. The author deduces Theorem A as a consequence of the 
Theorem B on cocycle superrigidity. 

Let the action of $\Gamma$ on $X$ be a free ergodic measure-preserving profinite action (i.e., an inverse 
limit of actions $\Gamma$ on $ X_n$ with $X_n$ finite) of a countable property $(T)$ group 
$\Gamma $ (more generally, of a group $\Gamma$ which admits an infinite normal subgroup $\Gamma_0$ such 
that the inclusion $\Gamma_0 \subset \Gamma$ has relative property $(T)$ and ${\Gamma} / {\Gamma_0}$  is finitely 
generated) on a standard probability space $X$. The author prove that if $\omega:\Gamma \times X \to \Lambda$ 
is a measurable cocycle with values in a countable group $\Lambda$, then $\omega$ is a cohomologous to a cocycle 
$\omega^{'}$ which factors through the map $\Gamma \times X \to \Gamma \times X_n$, for some $n$. As a corollary, 
he shows that any free ergodic measure-preserving action $\Lambda$ on $Y$ comes from a (virtual) 
conjugancy of actions. 

\section{Equidistribution for orbits of nonabelian semigroups on the torus}
Furstenberg \cite{Fu:DES} and Berent \cite{Be:MS} have investigated the action of abelian semigroups on the torus 
${\bf T}^d$ for $d=1$ and $d>1$ respectively. Their results answer problems raising by 
H. Furstenberg \cite{Fu:TI}.
Authors of the paper \cite{BFLM:SM}  extend to the noncommutative case some results of Furstenberg and Berent

\section{Zeta functions from spaces and dynamical systems}
Recall that Dedekind  has defined zeta function for polynomials over prime finite field. The zeta function is trivial and 
equal to $\frac{1}{1 - pz}$. However, combining the zeta function with           Chebyshev-Mobius  inversion formula we 
obtain the number of monic irreducible over ${\bf F}_p$ polynomials of natural degree $m$.   Riemann and Dedekind zeta 
functions are first examples of motivic zeta functions.
The authors of the paper \cite{FH:ZF} investigate sufficient conditions for (i) the existence of 
trace formulae for the Reidemeister number of a group endomorphism; 
(ii) the rationality of the Reidemeister zeta function and the 
convergence of the Nielsen zeta function; (iii) the equality of 
Reidemeister torsion of a group endomorphism to a special value of 
the Reidemeister zeta.     
   This interesting survey\cite{FH:ZF}  includes recent results on trace formulae, 
rationality and convergence of zeta functions and relations between 
special values of zeta functions and some simply homotopy invariants. 
The general setting of the paper \cite{MT:BZF}  is braided zeta functions in $q$-deformed geometry. In the framework 
authors define a zeta function for any rigid object in a ribbon braided category. In the ribbon case they define braided 
Hilbert series for objects in an Abelian braided category.
We will present some other types of zeta-functions.

\section{Dynamical Systems from Arithmetic Surfaces}

\subsection{Sato-Tate case}            
Let  $y^2 = f(x),f(x) = x^3 + c x  + d$   be a cubic polynomial in prime finite field ${\bf F}_{p}$.  For the number 
$\#C_p$ of points of the curve
$C: y^2 = f(x)$ in ${\bf F}_{p}$ the  well known formula 
$$ \#C_{p} = \sum_{x=0}^{p-1} \left( 1 + \left( \frac{f(x)}{p}\right) \right), $$
take place, where 
$\left( \frac{f(x_0)}{p}\right)$
is the Legendre symbol with a numerator which is equal to the value of the polynomial $f(x_0)$
in point $x_0 \in {\bf F}_{p}.$  It is ease to see that  $ \#C_{p}  = p - a_{p},$ where
$$a_{p} = - \sum_{x=0}^{p-1}  \left( \frac{f(x)}{p}\right)$$
If $C$ is the elliptic curve $ Å,$ then the number of points $\#C({\bf F}_{p})$  of the projective model
of the curve $Å$  in ${\bf F}_{p}$ is represented by the formula   $\#E_{p}  = 1 + p - a_p ,$
where $a_p = 2\sqrt p \cos\varphi_{p},$ If $C$ is not the elliptic curve, then the value $a_p$ 
is equal $1,\; -1$ or $0$ and ease to compute.
In both cases compute: 
$ \varphi_{p} = \arccos(a_p/2\sqrt p)$ and reduce it to the interval 
$[0, \pi].$ 

Let $E$ be an elliptic curve over rational numbers ${\bf Q}$ which does not admit complex multiplication.
Sato and Tate \cite{Birch} have given computational and theoretical evidences suggesting the distribution of
angles $ \varphi_{p}$.

Recently  L. Clozel, M. Harris, N. Shepherd-Barron, R. Taylor and their colleagues have proved the Sato-Tate
conjecture for all elliptic curves $E$ over ${\bf Q}$ (and over some its extensions) satisfying the mild condition
of having multiplicative reduction at some prime. 

Langlands conjectured that some symmetric power $L-$functions extend to
an entire function and coincide with certain automorphic $L-$functions. \\

{\bf Theorem} (Clozel, Harris, Shepherd-Barron, Taylor). Suppose $E$ is an elliptic
curve over $Q$ with non-integral $j-$invariant. Then for all $n > 0, L(s, E, Sym^{n})$
extends to a meromorphic function which is holomorphic and non-vanishing for
$Re(s) \ge 1 + n / 2 .$  \\

These conditions suffice to prove the Sato-Tate conjecture.

Theoretical considerations give

{\bf Proposition EC.} It is possible the arithmetic modeling of the Brownian motion by quantity $a_{p}$.

\subsection{ Kloosterman sums} 

 Let
  $$ T_{p}(c,d) = \sum_{x=1}^{p-1} e^{2{\pi}i(\frac{cx + \frac{d}{x}}{p})} $$

 $$ 1 \leq c, d \leq p - 1; \; x, c, d \in {\bf F}_{p}^{\ast} $$
 be a Kloosterman sum. \\
 By A. Weil estimate

  $$ T_{p}(c,d) = 2 \, \sqrt{p} \cos \theta_{p}(c,d) $$

  There are possible two distributions  of angles $\theta_{p}(c,d)$
on semiinterval $ [0, \pi ):$ \\
\\
a) {\it $p$ is fixed and $c$ and $d $ varies over ${\bf F}_{p}^{\ast};$
what is the distribution of angles  $ \theta_{p}(c,d)$ as
$p \rightarrow \infty \; ;$ }  \\
 \\
 b){\it $ \; c$ and $d$ are fixed and $p$ varies over all primes not
dividing $c$ and} $d.$ \\
\\
For the case $a)$ N. Katz~\cite{K:GS} and A. Adolphson~\cite{A:E} proved
that $\theta$ are distributed on $[0,\pi)$ with density
$\frac{2}{\pi} \sin^{2} t.$ \\
Let
$$ cd \not\equiv 0 \bmod p,  \; T_{p}(c,d) = \sum_{x=1}^{p-1} e^{2{\pi}i(\frac{cx + \frac{d}{x}}{p})}  $$
the Kloosterman sum. By A. Weil,
$T_{p}(c,d) = 2 \, \sqrt{p} \cos \theta_{p}(c,d)  .$
Compute $T_{p}, \cos \theta_{p}, \theta_{p}$ and reduce  $\theta_{p}$
to the interval $[0, \pi].$
Experiments demonstrate random behavior of angles of Kloosterman sums.

Theoretical considerations give 

{\bf Proposition KS.} It is possible the arithmetic modeling of the Brownian motion by Kloosterman sums.\\

{\Large \bf Conclusions}

We have presented a review of new results on actions and transformations of (quantum) groups and semigroups on 
(quantum) spaces, have presented dynamical systems  and   zeta functions arising from these spaces,  actions 
and transformations, discussed their stochastic properties.

\end{document}